\documentclass[12pt]{amsart}
\usepackage{amssymb,amsmath,amsfonts,latexsym,setspace}
\usepackage{bm}
\newcommand{\newsection}[1]{\setcounter{equation}{0} \section{#1}}
\newcommand{\bea}{\begin{eqnarray}}
\newcommand{\eea}{\end{eqnarray}}

\newcommand{\clb}{\mathcal{B}}

\newcommand{\cld}{\mathcal{D}}
\newcommand{\cle}{\mathcal{E}}
\newcommand{\clf}{\mathcal{F}}

\newcommand{\clh}{\mathcal{H}}

\newcommand{\cll}{\mathcal{L}}
\newcommand{\clm}{\mathcal{M}}

\newcommand{\clq}{\mathcal{Q}}

\newcommand{\cls}{\mathcal{S}}

\newcommand{\clw}{\mathcal{W}}

\newcommand{\raro}{\rightarrow}

\def \qed {\hfill \vrule height6pt width 6pt depth 0pt}
\def\textmatrix#1&#2\\#3&#4\\{\bigl({#1 \atop #3}\ {#2 \atop #4}\bigr)}
\def\dispmatrix#1&#2\\#3&#4\\{\left({#1 \atop #3}\ {#2 \atop #4}\right)}
\newcommand{\be}{\begin{equation}}
\newcommand{\ee}{\end{equation}}
\newcommand{\ben}{\begin{eqnarray*}}
\newcommand{\een}{\end{eqnarray*}}

\newcommand{\NI}{\noindent}

\newcommand{\bi}{\begin{itemize}}
\newcommand{\ei}{\end{itemize}}

\newtheorem{Theorem}{\sc Theorem}[section]
\newtheorem{Lemma}[Theorem]{\sc Lemma}
\newtheorem{Proposition}[Theorem]{\sc Proposition}
\newtheorem{Corollary}[Theorem]{\sc Corollary}
\newtheorem{Definition}[Theorem]{\sc Definition}
\newtheorem{Example}[Theorem]{\sc Example}
\newtheorem{Remark}[Theorem]{\sc Remark}
\newtheorem{Note}[Theorem]{\sc Note}
\newtheorem{Question}{\sc Question}
\newtheorem{ass}[Theorem]{\sc Assumption}
\newcommand{\bt}{\begin{Theorem}}
\def\beginlem{\begin{Lemma}}
\def\beginprop{\begin{Proposition}}
\def\begincor{\begin{Corollary}}
\def\begindef{\begin{Definition}}
\def\beginexamp{\begin{Example}}
\def\beginrem{\begin{Remark}}
\def\beginq{\begin{Question}}
\def\beginass{\begin{ass}}
\def\beginnote{\begin{Note}}
\newcommand{\et}{\end{Theorem}}
\def\endlem{\end{Lemma}}
\def\endprop{\end{Proposition}}
\def\endcor{\end{Corollary}}
\def\enddef{\end{Definition}}
\def\endexamp{\end{Example}}
\def\endrem{\end{Remark}}
\def\endq{\end{Question}}
\def\endass{\end{ass}}
\def\endnote{\end{Note}}

\begin{document}

\title[Invariant subspaces and reproducing kernel Hilbert spaces - I]{An Invariant Subspace Theorem and Invariant Subspaces of Analytic Reproducing Kernel Hilbert Spaces - I}

\author[Jaydeb Sarkar]{Jaydeb Sarkar}
\address{Indian Statistical Institute, Statistics and Mathematics Unit, 8th Mile, Mysore Road, Bangalore, 560059, India}
\email{jay@isibang.ac.in, jaydeb@gmail.com}

\subjclass[2010]{30H05, 46E22, 46M05, 46N99, 47A20, 47A45, 47B32,
47B38}


\keywords{Reproducing kernels, Hilbert modules, invariant subspaces,
weighted Bergman spaces, Hardy space}

\begin{abstract}
Let $T$ be a $C_{\cdot 0}$-contraction on a Hilbert space $\clh$ and
$\cls$ be a non-trivial closed subspace of $\clh$. We prove that
$\cls$ is a $T$-invariant subspace of $\clh$ if and only if there
exists a Hilbert space $\cld$ and a partially isometric operator
$\Pi : H^2_{\cld}(\mathbb{D}) \raro \clh$ such that $\Pi M_z = T
\Pi$ and that $\cls = \mbox{ran~} \Pi$, or equivalently, \[P_{\cls}
= \Pi \Pi^*.\]As an application we completely classify the
shift-invariant subspaces of analytic reproducing kernel Hilbert
spaces over the unit disc. Our results also includes the case of
weighted Bergman spaces over the unit disk.
\end{abstract}

\maketitle

\newsection{Introduction}

One of the most famous open problems in operator theory and function
theory is the so-called invariant subspace problem: Let $T$ be a
bounded linear operator on a Hilbert space $\clh$. Does there exist
a proper non-trivial closed subspace $\cls$ of $\clh$ such that $T
\cls \subseteq \cls$?

A paradigm is the well-known fact, due to Beurling, Lax and Halmos
(see \cite{B}, \cite{L}, \cite{H} and \cite{NF}), that any
shift-invariant subspace of $H^2_{\cle}(\mathbb{D})$ is given by an
isometric, or partially isometric, image of a vector-valued Hardy
space. Moreover, the isometry, or the partial isometry, can be
realized as an operator-valued bounded holomorphic function on
$\mathbb{D}$. More precisely, let $\cls$ be a non-trivial closed
subspaces of $H^2_{\cle}(\mathbb{D})$. Then $\cls$ is
shift-invariant if and only if there exists a Hilbert space $\cle_*$
such that $\cls$ is the range of an isometric, or partially
isometric operator from $H^2_{\cle_*}(\mathbb{D})$ to
$H^2_{\cle}(\mathbb{D})$ which intertwine the shift operators (see
\cite{RR}). Here $\cle$ is a separable Hilbert space and
$H^2_{\cle}(\mathbb{D})$ denote the $\cle$-valued Hardy space over
the unit disc $\mathbb{D} = \{z \in \mathbb{C} : |z| < 1\}$ (see
\cite{NF}, \cite{RR}).

In this paper we extend the Beurling-Lax-Halmos theorem for shift
invariant subspaces of vector-valued Hardy spaces to the context of
invariant subspaces of arbitrary $C_{\cdot 0}$-contractions. Recall
that a contraction $T$ on $\clh$ (that is, $\|Tf\| \leq \|f\|$ for
all $f \in \clh$) is said to be a $C_{\cdot 0}$-contraction if
$T^{*m} \raro 0$ as $m \raro \infty$ in the strong operator
topology. One of our main results, Theorem \ref{inv}, state that:
Let $\cls$ be a non-trivial closed subspace of a Hilbert space
$\clh$ and $T \in \clb(\clh)$ be a $C_{\cdot 0}$-contraction. Then
$\cls$ is a $T$-invariant subspace of $\clh$ if and only if there
exists a Hilbert space $\cle$ and a partial isometry $\Pi :
H^2_{\cle}(\mathbb{D}) \raro \clh$ such that $\Pi M_z = T \Pi$ and
that $\cls = \Pi H^2_{\cle}(\mathbb{D})$. This theorem will be
proven in Section 2.

In Section 3 we specialize to the case of reproducing kernel Hilbert
spaces, in which $T = M_z \otimes I_{\cle_*}$ and $\clh = \clh_K
\otimes \cle_*$. Here $\clh_K$ is an analytic Hilbert space (see
Definition \ref{ahs}) and $\cle_*$ is a coefficient space. In
Theorem \ref{MT} we show that a non-trivial closed subspace $\cls$
of $\clh_K \otimes \cle_*$ is $M_z \otimes I_{\cle_*}$-invariant if
and only if there exists a Hilbert space $\cle$ and a partially
isometric multiplier $\Theta \in \clm(H^2(\mathbb{D}) \otimes \cle,
\clh_{K} \otimes \cle_*)$ such that \[\cls = \Theta
H^2_{\cle}(\mathbb{D}).\]This classification extends the results of
Olofsson, Ball and Bolotnikov in \cite{O}, \cite{BV1} and \cite{BV2}
on the shift-invariant subspaces of vector-valued weighted Bergman
spaces with integer weights to that of vector-valued analytic
Hilbert spaces.

Our approach has two main ingredients: the Sz.-Nagy and Foias
dilation theory \cite{NF} and Hilbert module approach to operator
theory \cite{DP}. However, to avoid technical complications we speak
here just of operators on Hilbert spaces instead of Hilbert modules
over function algebras.

Finally, it is worth mentioning that in the study of invariant
subspaces of bounded linear operators on Hilbert spaces we lose no
generality if we restrict our attention to the class of $C_{\cdot
0}$-contractions.

\NI\textsf{Notations:} (1) All Hilbert spaces considered in this
paper are separable and over $\mathbb{C}$. We denote the set of
natural numbers including zero by $\mathbb{N}$. (2) Let $\clh$ be a
Hilbert space and $\cls$ be a closed subspace of $\clh$. The
orthogonal projection of $\clh$ onto $\cls$ is denoted by
$P_{\cls}$. (3) Let $\clh_1, \clh_2$ and $\clh$ be Hilbert spaces.
We denote by $\clb(\clh_1, \clh_2)$ the set of all bounded linear
operators from $\clh_1$ to $\clh_2$ and $\clb(\clh) = \clb(\clh,
\clh)$. (4) Two operators $T_1 \in \clb(\clh_1)$ and $T_2 \in
\clb(\clh_2)$ are said to be unitarily equivalent, denoted by $T_1
\cong T_2$, if there exists a unitary operator $U \in \clb(\clh_1,
\clh_2)$ such that $U T_1 = T_2 U$. (5) Let $\cle$ be a Hilbert
space. We will often identify $H^2_{\cle}(\mathbb{D})$ with
$H^2(\mathbb{D}) \otimes \cle$ and $M_z \in
\clb(H^2_{\cle}(\mathbb{D}))$ with $M_z \otimes I_{\cle} \in
\clb(H^2(\mathbb{D}) \otimes \cle)$, via the unitary $U \in
\clb(H^2_{\cle}(\mathbb{D}), H^2(\mathbb{D}) \otimes \cle)$ where
$U(z^m \eta) = z^m \otimes \eta$ and $m \in \mathbb{N}$ and $\eta
\in \cle$.

\section{An Invariant Subspace theorem}
In this section we will present a generalization of the
Beurling-Lax-Halmos theorem to the class of $C_{\cdot
0}$-contractions on Hilbert spaces.

Let $T$ be a $C_{\cdot 0}$-contraction on a Hilbert space $\clh$. A
fundamental theorem of Sz.-Nagy and Foias says that \[T \cong
P_{\clq} M_z|_{\clq},\]where $\clq$ is a $M_z^*$-invariant subspace
of $H^2_{\cld}(\mathbb{D})$ for some coefficient Hilbert space
$\cld$. In the following, we state and prove a variant of this fact
which is adapted to our present purposes.

\begin{Theorem}\label{PIH}
Let $T$ be a $C_{\cdot 0}$-contraction on a Hilbert space $\clh$.
Then there exists a coefficient Hilbert space $\cld$ and a
co-isometry $\Pi_{T} : H^2_{\cld}(\mathbb{D}) \raro \clh$ such that
$T \Pi_T = \Pi_T M_z$.
\end{Theorem}
\NI\textsf{Proof.} Let $D = (I_{\clh} - T T^*)^{\frac{1}{2}}$ and
$\cld = \overline{\mbox{ran}} D$. Since $\|z T^*\|_{\clb(\clh)} =
|z| \|T^*\|_{\clb(\clh)} < 1$, the inverse of $I_{\clh} - z T^*$
exists in $\clb(\clh)$ for all $z \in \mathbb{D}$.

\NI Define $L_T : \clh \raro H^2_{\cld}(\mathbb{D})$ by
\[(L_T h)(z):= D (I_{\clh} - z
T^*)^{-1} h = \sum_{m =0}^\infty ( D T^{*m} h) z^m. \quad \quad (h
\in \clh, \, z \in \mathbb{D})\] Now we compute
\[\begin{split} \|L_T h\|^2 & = \|\sum_{m = 0}^\infty ( D T^{*{m}} h) z^{m}\|^2 = \sum_{m = 0}^\infty \|D
T^{*m} h\|^2 = \sum_{m=0}^\infty \langle T^{m} D^2 T^{*m} h, h
\rangle
\\ & = \sum_{m=0}^\infty \langle T^m (I_{\clh} - T T^*) T^{*m}
h, h \rangle = \sum_{m=0}^\infty (\|T^{*m} h\|^2 - \|T^{*(m+1)}
h\|^2) \\& = \|h\|^2 - \lim_{m \raro \infty} \|T^{*m}
h\|^2,\end{split}\]where the last equality follows from the fact
that the sum is a telescoping series. This and the fact that
$\lim_{m \raro \infty} T^{*m} = 0$, in the strong operator topology,
implies that
\[\|L_T h\| = \|h\|. \quad \quad \quad (h \in \clh)\]Thus $L_T$ is
an isometry and $\Pi_T : H^2_{\cld}(\mathbb{D}) \raro \clh$ defined
by $\Pi_T = L_T^*$ is a co-isometry. Finally, for all $h \in \clh,\,
\eta \in \cld$ and $m \in \mathbb{N}$ we have
\[\begin{split}\langle \Pi_T (z^m \eta), h \rangle_\clh & = \langle
z^m \eta, D (I_{\clh} - z T^*)^{-1}
h\rangle_{H^2_{\cld}(\mathbb{D})} = \langle z^m \eta, \sum_{l \in
\mathbb{N}} (D T^{*l} h) z^{l} \rangle_{H^2_{\cld}(\mathbb{D})} \\
&= \langle \eta, D T^{*m} h\rangle_{\clh} = \langle T^{m} D \eta, h
\rangle_{\clh},\end{split}\]that is, \[\Pi_T (z^m \eta) = T^m D
\eta.\]Therefore, for each $m \in \mathbb{N}$ and $\eta \in \cld$ we
have \[\Pi_T M_z (z^m \eta) = \Pi_T (z^{m+1} \eta) = T^{m+1} D \eta
= T (T^m D \eta) = T \Pi_T (z^m \eta).\]Since $\{z^m \eta : m \in
\mathbb{N}, \eta \in \cld\}$ is total in $H^2_{\cld}(\mathbb{D})$,
it follows that
\[\Pi_{T} M_z = T \Pi_T.\]The proof is now complete. \qed

Now we present the main theorem of this section.

\begin{Theorem}\label{inv}
Let $T \in \clb(\clh)$ be a $C_{\cdot 0}$-contraction and $\cls$ be
a non-trivial closed subspace of $\clh$. Then $\cls$ is a
$T$-invariant subspace of $\clh$ if and only if there exists a
Hilbert space $\cld$ and a partially isometric operator $\Pi :
H^2_{\cld}(\mathbb{D}) \raro \clh$ such that \[\Pi M_z = T \Pi,\]and
that
\[\cls = \mbox{ran~} \Pi.\]
\end{Theorem}

\NI\textsf{Proof.} If $\cls$ is a $T$-invariant subspace of $\clh$
then \[(T|_{\cls})^{*m} = P_{\cls} T^{*m}|_{\cls} = P_{\cls} T^{*m},
\] shows that\[\|(T|_{\cls})^{*m} f \| = \|P_{\cls} T^{*m} f\| \leq
\|T^{*m} f\|,\]for all $f \in \cls$ and $m \in \mathbb{N}$. Thus
$T|_{\cls} \in \clb(\cls)$ is a $C_{\cdot 0}$-contraction. Now
Theorem \ref{PIH} implies that there exists a Hilbert space $\cld$
and a co-isometric map
\[\Pi_{T|_{\cls}} : H^2_{\cld}(\mathbb{D}) \raro \cls,\]such that
\[\Pi_{T|_{\cls}} M_z = T|_{\cls} \Pi_{T|_{\cls}}.\] Obviously
the inclusion map $i : \cls \raro \clh$ is an isometry and \[i
T|_{\cls} = T i.\]Let $\Pi$ be the bounded linear map from
$H^2_{\cld}(\mathbb{D})$ to $\clh$ defined by $\Pi = i
\Pi_{T|_{\cls}}$. Then
\[\Pi \Pi^* = (i \, \Pi_{T|_{\cls}}) (\Pi_{T|_{\cls}}^* \, i^*) = i i^* =
P_{\cls}.\]Therefore $\Pi$ is a partial isometry and $\mbox{ran}
\,\Pi = \cls$. Finally,
\[\Pi M_z = i \, \Pi_{T|_{\cls}} M_z = i T|_{\cls} \Pi_{T|_{\cls}} =
T i \Pi_{T|_{\cls}} = T \Pi.\] \NI This proves the necessary part.

\NI To prove the sufficient part it is enough to note that
$\mbox{ran} \Pi$ is a closed subspace of $\clh$ and $T \Pi = \Pi
M_z$ implies that $\mbox{ran} \Pi$ is a $T$-invariant subspace of
$\clh$. This completes the proof. \qed

The following corollary is a useful variation of the invariant
subspace theorem:
\begin{Corollary}\label{c1}
Let $T \in \clb(\clh)$ be a $C_{\cdot 0}$-contraction and $\cls$ be
a non-trivial closed subspace of $\clh$. Then $\cls$ is a
$T$-invariant subspace of $\clh$ if and only if there exists a
Hilbert space $\cld$ and a bounded linear operator $\Pi :
H^2_{\cld}(\mathbb{D}) \raro \clh$ such that $\Pi M_z = T \Pi$ and
\[P_{\cls} = \Pi \Pi^*.\]
\end{Corollary}

Before we go into the general theory of invariant subspaces of
reproducing kernel Hilbert spaces, let us consider the classical
Beurling-Lax-Halmos theorem as a simple corollary of Theorem
\ref{inv}.

\begin{Corollary}\label{BLH-D}
Let $\cls$ be a non-trivial closed subspace of the Hardy space
$H^2_{\cle}(\mathbb{D})$. Then $\cls$ is $M_z$-invariant if and only
if there exists a Hilbert space $\clf$ and a multiplier $\Theta \in
H^\infty_{\cll(\clf, \cle)}(\mathbb{D})$ such that $M_\Theta :
H^2_\clf(\mathbb{D}) \raro H^2_\cle(\mathbb{D})$ is partially
isometric and $\cls = \Theta H^2_{\clf}(\mathbb{D})$.
\end{Corollary}
\NI\textsf{Proof.} Let $\clf$ be a Hilbert space and $X :
H^2_{\clf}(\mathbb{D}) \raro H^2_{\cle}(\mathbb{D})$ be a bounded
linear map. It is easy to see that $X (M_z \otimes I_{\clf}) = (M_z
\otimes I_{\cle}) X$ if and only if $X = M_{\Theta}$, the
multiplication operator by $\Theta$, for some multiplier $\Theta \in
H^\infty_{\cll(\clf, \cle)}(\mathbb{D})$. Now the result follows
directly from Theorem \ref{inv}. \qed

Finally, note that Theorem \ref{inv} can be viewed as a
generalization of the Beurling-Lax-Halmos theorem for
shift-invariant subspaces of vector-valued Hardy spaces to invariant
subspaces of $C_{\cdot 0}$-contractions on Hilbert spaces.

\section{Analytic Hilbert spaces}

Let $\clh$ be a reproducing kernel Hilbert space of $\cle$-valued
holomorphic functions on $\mathbb{D}$ such that the multiplication
operator by the coordinate function, denoted by $M_z$, is bounded on
$\clh$ (cf. \cite{A}). A closed subspace $\cls$ of $\clh$ is said to
be \textit{shift-invariant} provided the product $z f \in \cls$
whenever $f \in \cls$.

The most general result on shift invariant subspaces has recently
been obtained by Olofsson, Ball and Bolotnikov in \cite{O},
\cite{BV1} and \cite{BV2}. Namely, for a given Hilbert space
$\cle_*$, a closed subspace $\cls$ of the weighted Bergman space
$L^2_{a, m}(\mathbb{D}) \otimes \cle_*$ ($m \geq 2$ and $m \in
\mathbb{N}$) is shift-invariant if and only if there exists a
Hilbert space $\cle$ and a function $\Theta : \mathbb{D} \raro
\clb(\cle, \cle_*)$ such that $M_{\Theta} : H^2(\mathbb{D}) \otimes
\cle \raro L^2_{a, m}(\mathbb{D}) \otimes \cle_*$ is a multiplier
(see definition below) and $\cls = \Theta H^2_{\cle}(\mathbb{D})$.
Recall that the weighted Bergman space $L^2_{a,
\alpha}(\mathbb{D})$, with $\alpha > 1$, is a reproducing kernel
Hilbert space corresponding to the kernel
\[k_{\alpha}(z, w) = \frac{1}{(1 - z \bar{w})^{\alpha}}. \quad \quad
(z, w \in \mathbb{D})\]

The purpose of this section is to extend the results of Olofsson,
Ball and Bolotnikov to a large class of reproducing kernel Hilbert
spaces. Our setting is very general and, as particular cases, we
obtain new and simple proof of the invariant subspace theorem for
vector-valued weighted Bergman spaces of Ball and Bolotnikov.

Let $K : \mathbb{D} \times \mathbb{D} \raro \mathbb{C}$ be a
positive definite function which is holomorphic in the first
variable, and anti-holomorphic in the second variable. We denote by
$\clh_K$ the reproducing kernel Hilbert space corresponding to the
kernel $K$.

\begin{Definition}\label{ahs} Let $\clh_K$ be a reproducing kernel Hilbert
space with $K$ as above. We say that $\clh_K$ is an analytic Hilbert
space if $M_z$ on $\clh_K$, defined by $M_z f = z f$ for all $f \in
\clh_K$, is a $C_{\cdot 0}$-contraction.
\end{Definition}

Let $\clh_{K_1}$ and $\clh_{K_2}$ be two analytic Hilbert spaces and
$\cle_1$ and $\cle_2$ be two Hilbert spaces. A map $\Theta :
\mathbb{D} \raro \clb(\cle_1, \cle_2)$ is said to be a
\textit{multiplier} from $\clh_{K_1} \otimes \cle_1$ to $\clh_{K_2}
\otimes \cle_2$ if
\[\Theta f \in \clh_{K_2} \otimes \cle_2. \quad \quad (f \in
\clh_{K_1} \otimes \cle_1)\]We denote the set of all multipliers
from $\clh_{K_1} \otimes \cle_1$ to $\clh_{K_2} \otimes \cle_2$ by
$\clm(\clh_{K_1} \otimes \cle_1, \clh_{K_2} \otimes \cle_2)$.

The following lemma, on a characterization of intertwining operators
between a vector-valued Hardy space and an analytic Hilbert space,
is well-known, which we prove for the sake of completeness.

\NI We will denote by $\mathbb{S}$ the Szego kernel on $\mathbb{D}$,
that is,
\[\mathbb{S}(z, w) = \frac{1}{(1 - z \bar{w})}. \quad \quad \quad
(z, w \in \mathbb{D})\]

\begin{Lemma}\label{lemma}
Let $\cle_1$ and $\cle_2$ be two Hilbert spaces and $\clh_K$ be an
analytic Hilbert space. Let $X \in \clb(H^2(\mathbb{D}) \otimes
\cle_1, \clh_{K} \otimes \cle_2)$. Then
\[X(M_z \otimes I_{\cle_1}) = (M_z \otimes I_{\cle_2}) X,\]if and
only if $X = M_{\Theta}$ for some $\Theta \in \clm(H^2(\mathbb{D})
\otimes \cle_1, \clh_{K} \otimes \cle_2)$.
\end{Lemma}
\NI\textsf{Proof.} Let $X \in \clb(H^2(\mathbb{D}) \otimes \cle_1,
\clh_{K} \otimes \cle_2)$ and $X(M_z \otimes I_{\cle_1}) = (M_z
\otimes I_{\cle_2}) X$. If $\zeta \in \cle_2$ and $w \in \mathbb{D}$
then
\[(M_z \otimes I_{\cle_1})^* [X^* (K(\cdot, w) \otimes
\zeta)] = X^* (M_z \otimes I_{\cle_2})^* (K(\cdot, w) \otimes \zeta)
= \bar{w} [X^* (K(\cdot, w) \otimes \zeta)],\]that is,
\[X^* (K(\cdot, w) \otimes \zeta) \in \mbox{ker} (M_z \otimes I_{\cle_1} -
wI_{H^2(\mathbb{D}) \otimes \cle_1})^*.\] This and the fact that
$\mbox{ker}(M_z - w I_{H^2(\mathbb{D})})^* = <\mathbb{S}(\cdot, w)>$
readily implies that
\[X^* (K(\cdot, w) \otimes \zeta) = \mathbb{S}(\cdot, w) \otimes
X(w) \zeta, \quad \quad \quad (w \in \mathbb{D}, \zeta \in
\cle_2)\]for some linear map $X(w) : \cle_2 \raro \cle_1$. Moreover,
\[\|X(w) \zeta\|_{\cle_1} = \frac{1}{\|\mathbb{S}(\cdot, w)\|_{H^2(\mathbb{D})}} \|X^* (K(\cdot, w) \otimes
\zeta)\|_{H^2(\mathbb{D}) \otimes \cle_1} \leq \frac{\|K(\cdot,
w)\|_{\clh_{K}}}{\|\mathbb{S}(\cdot, w)\|_{H^2(\mathbb{D})}} \|X\|
\|\zeta\|_{\cle_2},\]for all $w \in \mathbb{D}$ and $\zeta \in
\cle_2$. Therefore $X(w)$ is bounded and $\Theta(w) := X(w)^*  \in
\clb(\cle_1, \cle_2)$ for each $w \in \mathbb{D}$. Thus\[X^*
(K(\cdot, w) \otimes \zeta) = \mathbb{S}(\cdot, w) \otimes
\Theta(w)^* \zeta. \quad \quad \quad (w \in \mathbb{D}, \,\zeta \in
\cle_2)\] In order to prove that $\Theta(w)$ is holomorphic we
compute
\[\begin{split}\langle \Theta(w) \eta, \zeta \rangle_{\cle_2} & = \langle \eta, \Theta(w)^*
\zeta \rangle_{\cle_1} = \langle \mathbb{S}(\cdot, 0) \otimes \eta,
\mathbb{S}(\cdot, w) \otimes \Theta(w)^* \zeta
\rangle_{H^2(\mathbb{D}) \otimes \cle_1}
\\& = \langle X(\mathbb{S}(\cdot, 0) \otimes \eta), K(\cdot, w) \otimes
\zeta \rangle_{\clh_{K} \otimes \cle_2}. \quad \quad (\eta \in
\cle_1, \, \zeta \in \cle_2)\end{split}\]Since $w \mapsto K(\cdot,
w)$ is anti-holomorphic, we conclude that $w \mapsto \Theta(w)$ is
holomorphic. Hence $\Theta \in \clm(H^2(\mathbb{D}) \otimes \cle_1,
\clh_{K} \otimes \cle_2)$ and $X = M_{\Theta}$.

\NI Conversely, let $\Theta \in \clm(H^2(\mathbb{D}) \otimes \cle_1,
\clh_{K} \otimes \cle_2)$. For $f \in H^2(\mathbb{D}) \otimes
\cle_1$ and $w \in \mathbb{D}$ this implies that \[(z \Theta f)(w) =
w \Theta(w) f(w) = \Theta(w) w f(w) = (\Theta z f)(w).\]So
$M_{\Theta}$ intertwines the multiplication operators which
completes the proof. \qed

Now we are ready for the main theorem of this section.

\begin{Theorem}\label{MT}
Let $\clh_K$ be an analytic Hilbert space and $\cle_*$ be a Hilbert
space. Let $\cls$ be a non-trivial closed subspace of $\clh_K
\otimes \cle_*$. Then $\cls$ is $(M_z \otimes I_{\cle_*})$-invariant
subspace of $\clh_K \otimes \cle_*$ if and only if there exists a
Hilbert space $\cle$ and a partially isometric multiplier $\Theta
\in \clm(H^2(\mathbb{D}) \otimes \cle, \clh_{K} \otimes \cle_*)$
such that
\[\cls = \Theta (H^2(\mathbb{D}) \otimes \cle).\]
\end{Theorem} \NI\textsf{Proof.} By Theorem \ref{inv}, there exists a partial isometry
$\Pi : H^2(\mathbb{D}) \otimes {\cle} \raro \clh_K \otimes \cle_*$
such that $\Pi (M_z \otimes I_{\cle}) = (M_z \otimes I_{\cle_*})
\Pi$. Consequently, by Lemma \ref{lemma} we have that $\Pi =
M_{\Theta}$ for some $\Theta \in \clm(H^2(\mathbb{D}) \otimes
\cle_1, \clh_{K_2} \otimes \cle_2)$.

\NI The converse part is trivial. This completes the proof. \qed

In the present context, we restate Corollary \ref{c1} as follows:

\begin{Corollary}
Let $\clh_K$ be an analytic Hilbert space and $\cls$ be a
non-trivial closed subspace of $\clh_K \otimes \cle_*$ for some
coefficient Hilbert space $\cle_*$. Then $\cls$ is $(M_z \otimes
I_{\cle_*})$-invariant subspace of $\clh_K \otimes \cle_*$ if and
only if there exists a Hilbert space $\cle$ and a multiplier $\Theta
\in \clm(H^2(\mathbb{D}) \otimes \cle, \clh_{K} \otimes \cle_*)$
such that
\[P_{\cls} = M_{\Theta} M_{\Theta}^*.\]
\end{Corollary}

For each $\alpha > 1$, the weighted Bergman space $L^2_{a,
\alpha}(\mathbb{D})$ satisfies the conditions of Theorem \ref{MT}.
In particular, Theorem \ref{MT} includes the result by Ball and
Bolotnikov \cite{BV1} for weighted Bergman spaces with integer
weights as special cases.

\section{Concluding remarks}
A bounded linear operator $T$ on a Hilbert space $\clh$ is said to
have the \textit{wandering subspace property} if $\clh$ is generated
by the subspace $\clw_T : = \clh \ominus T \clh$, that is,
\[\clh = [\clw_T] = \overline{\mbox{span}}\{T^m \clw_T : m \in
\mathbb{N}\}.\]In that case we say that $\clw_T$ is a wandering
subspace for $T$.

\NI An important consequence of the Beurling theorem \cite{B} states
that: given a non-trivial closed shift invariant subspace $\cls$ of
$H^2(\mathbb{D})$, the subspace $\clw_{M_z|_{\cls}} = \cls \ominus z
\cls$ is a wandering subspace for $M_z|_{\cls}$. The same conclusion
holds for the Bergman space \cite{ARS} and the weighted Bergman
space with weight $\alpha = 3$ \cite{Sh} but for $\alpha > 3$, the
issue is more subtle (see \cite{HP}, \cite{MR}). In particular,
partially isometric representations of $M_z$-invariant subspaces of
analytic Hilbert spaces seems to be a natural generalization of the
Beurling theorem concerning the shift invariant subspaces of the
Hardy space $H^2(\mathbb{D})$.

Finally, it is worth stressing that the main results of this paper
are closely related to the issue of factorizations of reproducing
kernels. As future work we plan to extend our approach to several
variables and address issues such as factorizations of kernel
functions and containment of shift-invariant subspaces of analytic
Hilbert spaces over general domains in $\mathbb{C}^n$.

\vspace{0.2in}

\NI \textsf{Acknowledgement:} We would like to thank Joseph A. Ball
and Scott McCullough for kindly and quickly answering some questions
concerning reproducing kernel Hilbert spaces. We also thank Bhaskar
Bagchi and Dan Timotin for numerous useful comments and suggestions
on a preliminary version of this paper.

\end{document}